\documentclass[10pt]{article}
\usepackage{amssymb}
\usepackage{amsmath}
\usepackage[dvips]{graphics}
\usepackage{epsfig}
\newtheorem{theorem}{Theorem}

\begin{document}
\date{}
\author{Aristides V. Doumas$^{1}$ and Vassilis G. Papanicolaou$^{2}$\footnote{
Department of Mathematics, National Technical University of Athens, Zografou Campus, 
157 80 Athens, GREECE,
$^{1}$aris.doumas@hotmail.com \quad $^{2}$papanico@math.ntua.gr}}
\title{Uniform versus Zipf distribution in a mixing collection process}
\maketitle
\begin{abstract}
We consider the following variant of the classic collector's problem: The family of coupon probabilities is the mixing of two subfamilies one of which is the \textit{uniform} family, while the other belongs to the well known \textit{Zipf family}. We obtain asymptotics for the expectation, the second rising moment, and the variance of the random variable $T_N$, namely the number of trials needed for all the $N$ types of coupons to be collected (at least once, with replacement)
as $N \rightarrow \infty$. It is interesting that the effect of the uniform subcollection on the asymptotics of the expectation of $T_N$ (at least up to the sixth term) appears only in the leading factor of the expectation of $T_N$. The limiting distribution of $T_N$ is derived as well. These results answer a question placed in a recent work of ours [\textit{Electron. J. Probab.} \textbf{18} (2012) 1--15]. %motivated by problems coming from computer science and ecology.%
\end{abstract}

\textbf{2010 AMS Mathematics Classification.} 60F05; 60F99.

\section{Introduction and motivation}
The \textbf{``coupon collector's problem'' (CCP)} pertains to a population whose members are of $N$ different \emph{types}. For $1 \leq j \leq N$ we denote by $p_j$ the probability that a member of the population is of type $j$, where $p_j > 0$ and $\sum_{j=1}^{N}p_{j}=1$. We refer to the $p_j$'s as the \textit{coupon probabilities}. The members of the population are sampled independently \textit{with replacement} (alternatively, the polulation is assumed very large) and their types are recorded. Naturally, one quantity of interest is the number of trials $T_N$ needed until all $N$ types are detected (at least once). CCP belongs to the family of the so-called urn problems and it has been studied extensively; see, e.g., \cite{DPSETS} and the references therein. Moreover, due to its applications in several areas of science, new variants keep arising.

Let $\beta :=\{b_{j}\}_{j=1}^{\infty }$ be a sequence of strictly positive numbers. Then, for each integer $N \geq 1$
one can create a probability measure
$\pi _N =\{p_1,...,p_N\}$ on the set of types $\{1,...,N\}$ by taking
\begin{equation}
p_j = \frac{b_j}{B_N},
\qquad \text{where}\quad
B_N = \sum_{j=1}^N b_j.
\label{0}
\end{equation}
In a recent work (see \cite{DPM}) the authors asked what happens in the average when the sequence $\beta$ is the ``union'' of two subsequences one of which is constant (this corresponds to a uniform subcollection of coupons), while the other obeys some rather general law, in particular the law of the well-known \textit{Zipf} family.\footnote{In \cite{DPM} the authors also asked the same question when the family of coupon probabilities is the ``mixing'' of two constant subsequences. For an answer see \cite{DM}.} \textit{Zipf}, this surprising law of nature, arises in many areas of science, such as computer science, physics, biology, earth and planetary sciences, economics and finance, as well as linguistics, demography, and the social sciences (see, e.g., the highly
cited article \cite{NEWMAN} of Mark Newman, where he reviewed some of the empirical evidence for the existence of
power-law forms, and the recent work \cite{LEN} of Locey and Lennon on the applications of power-laws in biology).

In this paper we bring an answer to the above question by deriving the asymptotics of the expectation and of the second moment (up to the fifth and sixth term respectively) as $N \to \infty$, as well as the limit distribution of $T_N$ (under the apropriate normalization). Let
\begin{equation}
b_{2j-1} = 1
\qquad \text{and} \qquad
b_{2j} =  a_j,
\qquad\qquad
j = 1, 2, \dots,
\label{sb1}
\end{equation}
where $\{a_j\}_{j=1}^{\infty } =: \alpha$ is a sequence of strictly positive numbers of the form
\begin{equation}
a_{j}=\frac{1}{j^{p}},\,\,\, p>0.
\label{aj}
\end{equation}
The case where $p=1$ corresponds to the \textit{standard Zipf} distribution. For general positive values of $p$ we have the so--called \textit{generalized Zipf} subfamily of coupons. 

Testing uniform and the standard Zipf distribution is not a new idea. We refer the reader to the highly cited articles \cite{PEER} on the search  and replication in unstructured peer-to-peer networks, and \cite{YCSB} on the benchmarking cloud serving systems with the Yahoo! Cloud
Serving Benchmark (YCSB) framework. However, in this paper we consider the problem of the coexistence of uniform and generalized \textit{Zipf} distributions in the \textit{same} model. The question about the effect of the uniform--Zipf distribution on the average of the random variable $T_{N}$ arises naturally. As we will see, the uniform subcollection acts on the asymptotics of the expectation of $T_N$ only in the leading factor (at least up to the fifth term of its asymptotic expansion). The same argument holds for the second rising moment of $T_N$ up to the sixth term. In comparison with the classic version of the problem (when all coupons are uniformly distributed), or with the case where all coupons are Zipf distributed, the effect of the uniform subcollection (in the mixing case studied here) causes a significant increment in the number of trials needed for a complete set of coupons. This argument will be illustrated via an example at the end of the paper.

\section{Main results}

It is well known (see, e.g., \cite{F-G-T}) that the expectation of $T_N$ can be expressed as
\begin{equation}
E\left[ T_N \right] = \int_0^{\infty} \left[1-\prod_{j=1}^{N}\left(1 - e^{-p_j t} \right)\right] dt
= \int_{0}^{1}\left[1-\prod_{j=1}^N \bigg(1-x^{p_{j}}\bigg)\right]\frac{dx}{x}.
\label{mv}
\end{equation}
From now on we assume that $N$ is even and for convenience we set
\begin{equation}
N := 2M.
\label{mv1}
\end{equation}
By substituting $t=-B_N\ln y$ and thanks to the binomial theorem,
formula (\ref{mv}) (in view of \eqref{sb1}) yields
\begin{equation}
E\left[ T_N \right] = B_N \int_0^{1} \left[1-\prod_{j=1}^{M}\left(1 - y^{a_j } \right)-\sum_{k=1}^{M}\binom {M}{k}(-1)^{k}y^{k}\prod_{j=1}^{M}\left(1 - y^{a_j } \right)\right] \frac{dy}{y}.\label{999}
\end{equation}
Notice that from (\ref{0}) and (\ref{sb1}) we have
\begin{equation}
B_N = M + A_M,
\qquad \text{where }\;
A_M :=\sum_{j=1}^M a_j.
\label{QAZ}
\end{equation}
The study of the quantity $A_M$ of (\ref{QAZ}) is an external matter. In particular, one easily gets its full asymptotic expansion via the celebrated Euler--Maclaurin summation formula, as we will shortly see in the last step of the proof of our main theorem. Let $\tilde{T}_M$ be the number of trials needed for one to collect (with replacement) all $M$ different types of coupons when the coupon probabilities are
\begin{equation*}
q_j := \frac{a_j}{A_M},
\qquad
j = 1, \dots, M.
\end{equation*}
Then, \eqref{mv} implies
\begin{equation}
E\left[\tilde{T}_M\right] = A_{M} \int_0^{1} \left[1-\prod_{j=1}^{M}\left(1 - y^{a_j } \right)\right] \frac{dy}{y}.\label{r1}
\end{equation}
Thus, (\ref{999}) yields
\begin{equation}
E\left[ T_N \right] =B_{N}\left[ A_M^{-1} E\left[\tilde{T}_M\right]
-\sum_{k=1}^M \binom {M}{k}(-1)^{k}\int_0^{1}y^{k-1}\prod_{j=1}^M \left(1 - y^{a_j } \right)dy\right].
\label{P7}
\end{equation}
The main results of the paper are presented in the following
\begin{theorem}
Let the sequence $\beta = \{b_{j}\}_{j=1}^{\infty }$ be the ``union'' of two subsequences, as given by (\ref{sb1}), (\ref{aj}), one of which is constant (this corresponds to a uniform subcollection of coupons), while the other
belongs to the generalized Zipf family, namely $\alpha = \{a_j = 1/j^{p},\,\, p > 0\}$. Then,  as $N=2M \rightarrow \infty$ we have
\begin{align}
E\left[T_N \right]= M^{p+1} &\left[\ln M-\ln\left(\ln\frac{M}{p}\right)+\left(\gamma-\ln p\right)+\frac{\ln\left(\ln\frac{M}{p}\right)}{\ln{M}}\right.\nonumber\\
&\left.-\frac{1+\gamma+\frac{1}{p}}{\ln{M}}+ O\left(\frac{\ln\left(\ln{M}\right)}{\ln{M}}\right)^{2}\right],
\label{MEAN 2}
\end{align}
where $\gamma$ is, as usual, the Euler--Mascheroni constant. Regarding the second rising moment\footnote{under the notation $t^{(2)} = t (t+1)$
} and the variance of the r.v. $T_{N}$ we have
\begin{align}
E\left[T_{N}^{(2)}\right]=& M^{2p+2} \left[\ln^{2} M+2\left(\gamma-\ln p\right)\ln M-2\ln\left(\ln\frac{M}{p}\right)\ln M+\left(\ln\left(\ln\frac{M}{p}\right)\right)^{2}\right.\nonumber\\
&\left.+2\left (\ln p -\gamma+1\right)\ln\left(\ln\frac{M}{p}\right)+\left(\gamma^{2}+\frac{\pi^{2}}{6}-2\gamma-2-\frac{2}{p}+\ln^{2}p\right)
+ O\left(\frac{\ln\left(\ln{M}\right)}{\ln{M}}\right)^{2}\right],
\label{IR2}
\end{align}
\begin{equation}
V\left[T_N\right] \sim \frac{\pi^2}{6}\;M^{2p+2}.
\label{FINAL}
\end{equation}
Moreover, $T_N$ appropriately normalized converges in distribution to a standard Gumbel random variable. More precisely as $N\rightarrow \infty$
\begin{equation}
P\left\{\frac{T_N - m_N}{k_N} \leq y \right\}\longrightarrow \exp(e^{-y})
\qquad
\text{for all}\ \, y \in \mathbb{R},\,\,\, N=2M,
\label{N8b}
\end{equation}
where,
\begin{equation}
m_N = M^{p+1}\left[\ln \left(\frac{M}{p}\right) - \ln\left(\ln \left(\frac{M}{p}\right)\right)\right] \quad  \text{and} \quad k_N = M^{p+1},\label{bk}
\end{equation}
\end{theorem}
\textit{Proof of Theorem 1.} Starting from (\ref{P7}) (recall that $a_{j}=j^{-p}, \,p>0$), we focus on the quantities
\begin{equation}
\,\,\,\,\,\,\,\,\,\,\,\,\,\,\,\,\,\,\,\,W_k \left(M\right) :=
\int_0^1 y^{k-1}\prod_{j=1}^M \left(1 - y^{j^{-p} } \right) dy,\,\,\,\,\,k=1,2,\dots, M.
\label{W}
\end{equation}
If we set
\begin{equation}
F(x) := x^{p} \ln\left(\frac{x}{p}\right),
\label{F}
\end{equation}
then, in view of (\ref{aj}) and under the change the variables $y=e^{-sF(M)}$ formula (\ref{W}) becomes
\begin{equation}
W_k \left(M\right) = M^{p}\ln\left(\frac{M}{p}\right)\int_0^{\infty} e^{-k s M^{p} \ln\left(\frac{M}{p}\right)}\prod_{j=1}^M \left(1 - e^{-s\left(\frac{M}{j}\right)^{p} \ln\left(\frac{M}{p}\right)
 } \right) ds.
\label{990}
\end{equation}
The following result is important for our analysis:
\begin{align}
\int_{1}^{M}e^{-s\left(\frac{M}{x}\right)^{p} \ln\left(\frac{M}{p}\right)}dx=
\frac{1}{s} \left(\frac{M}{p}\right)^{1-s} \frac{1}{\ln\left(\frac{M}{p}\right)}
&-\left(1+\frac{1}{p}\right)\frac{1}{s^{2}} \left(\frac{M}{p}\right)^{1-s} \frac{1}{\ln^{2}\left(\frac{M}{p}\right)} \nonumber \\
&\times \left[1+O\left(\frac{1}{\ln M}\right)\right], \label{IIa}
\end{align}
uniformly in $s\in [s_{0},\infty)$, for any fixed $s_{0} > 0 $.

The proof is based on the method of integration by parts and is omitted. By (\ref{IIa}), the comparison of sums and integrals, and the Taylor expansion of the logarithm we get
\begin{equation}
\lim_{M}\sum_{j=1}^{M}\ln \Bigg( 1-e^{-s\left(\frac{M}{j}\right)^{p} \ln\left(\frac{M}{p}\right)
}\Bigg) = \left\{
\begin{array}{rc}
-\infty,&  \text{if } s<1 \\
0,&   \text{if } s\geq1,
\end{array}
\right.  \label{SL1A}
\end{equation}
Taking advantage of (\ref{SL1A}) and for any given $\varepsilon \in (0,1)$ we rewrite (\ref{990}) as
\begin{align}
W_{k}\left(M;\alpha\right)= M^{p}\ln\left(\frac{M}{p}\right)\Bigg(\,I_1 (M)+I_2 (M)+I_3 (M)\,\Bigg),
\label{b3a}
\end{align}
where
\begin{align}
I_1 (M):&= \int_0^{1-\varepsilon}\left[ \exp \left\{-ks M^{p}\ln\left(\frac{M}{p}\right)+\sum_{j=1}^M \ln \Bigg( 1-e^{-\left(\frac{M}{j}\right)^{p} s\ln\left(\frac{M}{p}\right)
}\Bigg) \right\} \right] ds,\label{I1}\\
I_{2}(M): &= \int_{1-\varepsilon}^{1}\left[ \exp \left\{-ksM^{p}\ln\left(\frac{M}{p}\right)+\sum_{j=1}^M \ln \Bigg( 1-e^{-\left(\frac{M}{j}\right)^{p} s\ln\left(\frac{M}{p}\right)
}\Bigg) \right\} \right] ds,\label{I2}\\
I_{3}(M): &=\int_{1}^{\infty }\left[ \exp \left\{-ksM^{p}\ln\left(\frac{M}{p}\right)+\sum_{j=1}^M \ln \Bigg( 1-e^{-\left(\frac{M}{j}\right)^{p} s\ln\left(\frac{M}{p}\right)
}\Bigg) \right\} \right] ds.\label{I3}
\end{align}
As we will see all the information we need comes from $I_{2}(M)$. Starting from (\ref{I3}) and using  (\ref{SL1A}) we get
\begin{equation*}
I_{3}(M)=\int_{1}^{\infty }e^{-ksM^{p}\ln\left(\frac{M}{p}\right)}
\left\{1-\int_{1}^{M}e^{-s\left(\frac{M}{x}\right)^{p}\ln\left(\frac{M}{p}\right)}dx \left[1+O\left(\int_{1}^{M}e^{-s\left(\frac{M}{x}\right)^{p}\ln\left(\frac{M}{p}\right)}dx\right)\right]\right\}ds.
\end{equation*}
By invoking (\ref{IIa}) and integrating by parts the above becomes
\begin{align}
I_{3}(M)=\frac{1}{k M^{p}\ln\left(\frac{M}{p}\right)}
e^{-k M^{p}\ln\left(\frac{M}{p}\right)}
-&\frac{1}{k M^{p}\ln^{2}\left(\frac{M}{p}\right)}
e^{-k M^{p}\ln\left(\frac{M}{p}\right)} \nonumber \\
&\times\left[1+O\left(\frac{1}{k M^{p}\ln M}
e^{-k M^{p}\ln M}\right)\right],\,\,k=1,2,\cdots,M.\label{10}
\end{align}
Our next task is $I_{2}(M)$ of (\ref{I2}). By applying the Taylor expansion of the logarithm and using the comparison of sums and integrals, as well as the result presented in formula (\ref{IIa}) (since $s$ in this case is \textit{strictly} positive), and finally, changing the variables as
\begin{equation*}
u=\frac{1}{\ln\left(\frac{M}{p}\right)}\left(\frac{M}{p}\right)^{1-s}
\end{equation*}
one arrives at
\begin{align}
I_{2}(M)=&\frac{1}{\ln\left(\frac{M}{p}\right)}
 e^{-k M^{p} \ln\left[\frac{1}{\ln\left(\frac{M}{p}\right)}e^{\ln\left(\frac{M}{p}\right)}\right]}\nonumber\\
&\int_{1/\ln\left(\frac{M}{p}\right)}^{\left(\frac{M}{p}\right)^{\epsilon}/\ln\left(\frac{M}{p}\right)}
e^{k M^{p}\ln u}
\exp\left(-\frac{u}{1-\frac{\ln u}{\ln\left(\frac{M}{p}\right)}-\frac{\ln\left(\ln\left(\frac{M}{p}\right)\right)}{\ln\left(\frac{M}{p}\right)}}\right.\nonumber\\
&\left.+\frac{\left(1+\frac{1}{p}\right)u}{\ln\left(\frac{M}{p}\right)\left[1-\frac{\ln u}{\ln\left(\frac{M}{p}\right)}-\frac{\ln\left(\ln\left(\frac{M}{p}\right)\right)}{\ln\left(\frac{M}{p}\right)}\right]^{2}}
\left[1+O\left(\frac{1}{\ln M}\right)\right]\right)\frac{du}{u}.
\label{189}
\end{align}
Since, for $\left|x\right|<1$,  $\left(1-x\right)^{-2}=\sum_{n=1}^{\infty}n x^{n-1}$, the integral appearing in (\ref{189}) yields
\begin{align*}
&\int_{1/\ln\left(\frac{M}{p}\right)}^{\left(\frac{M}{p}\right)^{\epsilon}/\ln\left(\frac{M}{p}\right)}
\frac{e^{k M^{p}\ln u -u}}{u}
\exp\left(-u\sum_{n=1}^{\infty}
\left( \frac{1}{\ln\left(\frac{M}{p}\right)} \ln\left[u\ln\left(\frac{M}{p}\right)\right]\right)^n \right)\\
&\times\exp\left(\left(1+\frac{1}{p}\right)\frac{1}{\ln\left(\frac{M}{p}\right)}\, u
\left[1 + O\left(\frac{1}{\ln\left(\frac{M}{p}\right)}\right)\right]
\sum_{n=1}^{\infty} n \left( \frac{1}{\ln\left(\frac{M}{p}\right)} \ln\left[u\ln\left(\frac{M}{p}\right)\right]\right)^{n-1}\right) du.
\end{align*}
In order to obtain the leading behavior of the integral above as $N = 2M \rightarrow \infty$ it suffices to work with the integral
\begin{equation*}
J(M) := \int_{1/\ln\left(\frac{M}{p}\right)}^{\left(\frac{M}{p}\right)^{\epsilon}/\ln\left(\frac{M}{p}\right)}
e^{k M^{p}\ln u -u}\,\,\frac{du}{u}.
\end{equation*}
Changing the variables as $u=M^{p}s$ and applying the Laplace method for integrals (see, e.g., \cite{B-O}) we arrive at
\begin{equation*}
J(M)\sim \frac{1}{k}\, \frac{\ln\left(\frac{M}{p}\right)}{\ln\left(M^{p}\ln\left(\frac{M}{p}\right)\right)}\,
\, e^{-\frac{1}{\ln\left(\frac{M}{p}\right)}} \,\,
e^{k M^{p}\ln\left(\frac{1}{M^{p}\ln\left(\frac{M}{p}\right)}\right)}, \,\,\,\,M\rightarrow \infty
\end{equation*}
and by invoking (\ref{189}) one gets
\begin{equation}
I_{2}(M)\sim \frac{1}{k}\, \frac{1}{\ln\left(M^{p}\ln\left(\frac{M}{p}\right)\right)}\,
\, e^{-\frac{1}{\ln\left(\frac{M}{p}\right)}} \,\,
e^{-k M^{p}\ln\left(\frac{M}{p}\right)}, \,\,\,\,M\rightarrow \infty. \label{000}
\end{equation}
From (\ref{10}) and (\ref{000}) one has that $I_3 (M)$ is negligible compared to $I_2 (M)$ as $M\rightarrow \infty$. Finally, for $I_1 (M)$ of (\ref{I1}) we have
\begin{align*}
I_1 (M)&< \int_0^{1-\varepsilon}\left[ \exp \left\{\sum_{j=1}^M \ln \Bigg( 1-e^{-\left(\frac{M}{j}\right)^{p} s\ln\left(\frac{M}{p}\right)
}\Bigg) \right\} \right] ds\\
&< \exp \left(-\sum_{j=1}^M e^{-\left(\frac{M}{j}\right)^{p} \left(1 - \varepsilon \right)} \right)<\exp \left(-\int_{1}^{M} e^{-\left(\frac{M}{x}\right)^{p}\left(1-\varepsilon \right)}dx \right).
\end{align*}
From (\ref{IIa}) and (\ref{000}) one has that $I_1 (M)$ is negligible compared to $I_2 (M)$ as $M\rightarrow \infty$ and, as we have seen, the same argument holds for $I_3 (M)$. Hence, from (\ref{b3a}) we get
\begin{equation}
W_{k}\left(M\right)\sim \frac{1}{k}\, \frac{M^{p}\ln\left(\frac{M}{p}\right)}{\ln\left(M^{p}\ln\left(\frac{M}{p}\right)\right)}\,
\, e^{-\frac{1}{\ln\left(\frac{M}{p}\right)}} \,\,
e^{-k M^{p}\ln\left(\frac{M}{p}\right)}, \,\,\,\,M\rightarrow \infty. \label{1b3a}
\end{equation}
To complete our analysis, and in view of (\ref{P7}), one must obtain the leading term of the quantity
\begin{equation*}
\sum_{k=1}^{M}\binom {M}{k}(-1)^{k}W_{k}\left(M\right).
\end{equation*}
It is not hard to check that
\begin{equation}
\sum_{k=1}^{M}\binom {M}{k}(-1)^{k}W_{k}\left(M\right)\sim
-\frac{M^{p}\ln\left(\frac{M}{p}\right)}{\ln\left(M^{p}\ln\left(\frac{M}{p}\right)\right)}\,
\, e^{-\frac{1}{\ln\left(\frac{M}{p}\right)}} \,\,
e^{- M^{p}\ln\left(\frac{M}{p}\right)}, \,\,\,\,M\rightarrow \infty.\label{888}
\end{equation}
Let us now return to (\ref{P7}) and the quantity $E\left[\tilde{T}_M\right]$. Under (\ref{aj}) the first five terms of the asymptotics of $E\left[\tilde{T}_M\right]$ (as $M\rightarrow \infty$) are known. In particular, (see \cite{DP} and \cite{DPSETS})
\begin{align}
E\left[\tilde{T}_M\right] = A_M M^{p} &\left[\ln M-\ln\left(\ln\frac{M}{p}\right)+\left(\gamma-\ln p\right)+\frac{\ln\left(\ln\frac{M}{p}\right)}{\ln\frac{M}{p}}\right.\nonumber\\
&\left.-\frac{1+\gamma+\frac{1}{p}}{\ln\frac{M}{p}}+ O\left(\frac{\ln\left(\ln{M}\right)}{\ln{M}}\right)^{2}\right].
\label{R2}
\end{align}
By invoking (\ref{888}) and (\ref{R2}) in (\ref{P7}) we have
\begin{align}
E\left[T_N \right]= \left(M+\sum_{j=1}^{M}\frac{1}{j^{p}}\right)  M^{p} &\left[\ln M-\ln\left(\ln\frac{M}{p}\right)+\left(\gamma-\ln p\right)+\frac{\ln\left(\ln\frac{M}{p}\right)}{\ln {M}}\right.\nonumber\\
&\left.-\frac{1+\gamma+\frac{1}{p}}{\ln{M}}+ O\left(\frac{\ln\left(\ln{M}\right)}{\ln{M}}\right)^{2}\right].
\label{MEAN}
\end{align}
\textit{Last step before the expectation}. To obtain the asymptotics of $E\left[T_N \right]$ one has to investigate the asymptotics of $A_{M}=\sum_{j=1}^{M}j^{-p}$. By the celebrated Euler--Maclaurin summation formula (see, e.g. \cite{B-O}) the full aymptotic expansion of $A_{M}$ is known (as $M\rightarrow \infty$). In particular, the leading term in the asymptotics of $A_{M}$ depends on the behaviour of the series $\sum_{j=1}^{\infty}1/j^{p}$. If $p>1$ we have
\begin{equation}
A_{M}\sim \zeta (p) \label{101}
\end{equation}
where $\zeta (p)$ denotes the Riemann zeta function, while for $0<p< 1$ we have
\begin{equation}
A_{M}\sim \int_{1}^{M} x^{-p}dx=\frac{M^{1-p}}{1-p}.\label{102}
\end{equation}
For $p=1$, namely the case of the \textit{standard Zipf distribution} we have
\begin{equation}
A_{M}\sim \ln M.\label{103}
\end{equation}
\underline{\textbf{Claim.}} The effect of the uniform subcollection on the asymptotics of the expectation of $T_N$ (at least up to the sixth term) appears \textit{only} in the leading factor of \eqref{MEAN}. To wit (see \eqref{MEAN}) it suffices to check that as $M \to \infty$
\begin{equation}
M  \left(\frac{\ln\left(\ln{M}\right)}{\ln{M}}\right)^{2}>> A_{M}\ln M.\label{104}
\end{equation}
The proof of \eqref{104} is immediate in all three cases given in ((\ref{101})--(\ref{103})). 
The result for the expectation of the r.v. $T_N$ now follows by invoking (\ref{104}) in (\ref{MEAN}).\\
\textit{Second moment, variance and distribution of $T_{N}$}. Mimicking the derivation of the asymptotics of $E[T_{N}]$ it is straightforward to get the asymptotics of the second rising moment of the random variable $T_{N}$. We have (see, e.g., \cite{DP})
\begin{align}
E\left[T_{N}^{(2)}\right]=-2\int_{0}^{1}\left[1-\prod_{j=1}^N \bigg(1-x^{p_{j}}\bigg)\right]\frac{\ln x}{x} \, dx,
\label{14}
\end{align}
where we have used the notation $t^{(2)} = t (t+1)$. Similarly to formula (\ref{999}) we get
\begin{equation}
E\left[T_{N}^{(2)}\right]=-2 B^{2}_N \int_0^{1} \left[1-\prod_{j=1}^{M}\left(1 - y^{a_j } \right)-\sum_{k=1}^{M}\binom {M}{k}(-1)^{k}y^{k}\prod_{j=1}^{M}\left(1 - y^{a_j } \right)\right] \frac{\ln y}{y}dy.\label{1999}
\end{equation}
Likewise, similarly to (\ref{P7}) one has
\begin{equation}
E\left[T_{N}^{(2)}\right]=B^{2}_N\left[ A_M^{-2} E\left[\tilde{T}_M^{2}\right]
+2\sum_{k=1}^M \binom {M}{k}(-1)^{k}Q_k \left(M\right)\right],
\label{PP7}
\end{equation}
where
\begin{equation}
\,\,\,\,\,\,\,\,\,\,\,\,\,\,\,\,\,\,\,\,Q_k \left(M\right) :=
\int_0^1  y^{k-1}\,\ln y\prod_{j=1}^M \left(1 - y^{j^{-p} } \right) dy,\,\,\,\,\,k=1,2,\dots, M,
\label{Q}
\end{equation}
and as in (\ref{r1})
\begin{equation}
E\left[\tilde{T}_M^{2}\right] = -2 A^{2}_{M} \int_0^{1} \left[1-\prod_{j=1}^{M}\left(1 - y^{a_j } \right)\right] \frac{\ln y}{y}dy.\label{1r1}
\end{equation}
Under formula (\ref{aj}) the first six terms of the asymptotics of $E\left[\tilde{T}_M^{2}\right]$ (as $M\rightarrow \infty$) are known (see \cite{DP} and \cite{DPSETS}). Finally, one arives at the desired result. Again, the effect of the uniform subcollection in the asymptotics of the second rising moment of the random variable $T_{N}$ appears only in the leading factor of the second rising moment of the random variable $T_{N}$.\\
\textbf{Observation.} It is straightforward for one to check that the same result holds for all the rising moments of the random variable $T_N$. Having (\ref{MEAN 2}) and (\ref{IR2}) it is easy to obtain leading asymptotics for the variance of $T_{N}$. Using the formula
\begin{equation*}
V\left[T_N\right]=E\left[T_{N}^{(2)}\right]-E[\,T_{N}\,]-E[\,T_{N}\,]^{2}
\end{equation*}
we get (\ref{FINAL}) as $N\rightarrow \infty$.
The previous results drive us to normalize $T_N$ as
\begin{equation*}
\frac{T_N - m_N}{k_N}
\end{equation*}
where, $m_N$ and $k_{N}$ are given in (\ref{bk}), and by a well known theorem (see, e.g., \cite{DPSETS}) one obtains the final result of Theorem 1 (i.e., the r.v. $T_{N}$ converges in distribution to a standard Gumbel r.v.). We remind the reader that in the classic version of the problem (namely, the case of one class of uniformly distributed coupons) the corresponding limiting theorem is due to P. Erd\H{o}s and A. R\'{e}nyi:
\begin{equation}
P\left\{\frac{T_N - N\ln N}{N} \leq y \right\}\longrightarrow \exp(e^{-y})
\qquad
\text{for all}\ \, y \in \mathbb{R},
\label{N8c}
\end{equation}
see \cite{E-R}, while for the case the coupon probabilities are distributed according to the \textit{Zipf law} we have the following theorem (see \cite{DP} and \cite{DPSETS})
\begin{equation}
P\left\{\frac{T_N - A_{N}\,N^{p}\left[\ln \left(\frac{N}{p}\right) - \ln\left(\ln \left(\frac{N}{p}\right)\right)\right]}{A_{N}\,N^p} \leq y \right\}\longrightarrow \exp(e^{-y})
\qquad
\text{for all}\ \, y \in \mathbb{R},
\label{N8d}
\end{equation}
To support the above limiting results let us consider the following\\
\textbf{Example.} 
Recall that the coupon probabilities $p_j$ satisfy (\ref{0})-(\ref{sb1}), where $j=1,2,\cdots, N$ and $N=2M$. Let us compute the minimum number of trials, so that with probability $0.90$ we get a complete set of all $N$ different types of coupons when $N=2M=100$ and $a_{j}=1/j,\;j=1,2,\cdots,M$.

We have $N=100,\; f(M)=M$. Hence, $b_{100}=50^{2}\left(\ln(50)-\ln(\ln(50))\right)=6,369.92,\;k_{100}=50^{2}$. Assume that the answer is $a$ trials. By (\ref{N8b}) we have 
\begin{eqnarray*}
P\left(T^{\text{mix}}_{100}\leq a\right)&=&P\left(\left(T^{\text{mix}}-6,369.22\right)/2500\leq \left(a -6,369.22\right)/2500\right)  \\
&\approx&\exp(-e^{-\lambda})=0.90, 
\end{eqnarray*} 
where \, $\lambda=\left(a -6,369.22\right)/2500.$ So that $\lambda=-\ln\left[-\ln\left(0.90\right)\right]=2.25037$. Thus, with probability $0.90$ one needs at least $\textbf{11,996}$ trials to collect all $100$ different types of coupons.  
\smallskip

Now, let us compare our results with the classic version of the problem when all the $N$ different coupons are distributed according to the \textit{standard Zipf law}. In this case we have from (\ref{N8d}): $N=100,\; f(N)=N,\, A_{N}=H_{100}=5.18738.$ Hence, $b_{100}=1,596.67,\;k_{100}=518.738$. Assume that the answer is $k$ trials. We have 
\begin{eqnarray*}
P\left(T^{\text{Zipf}}_{100}\leq k\right)&=&P\left(\left(T^{\text{Zipf}}_{100}-1,596.67\right)/518.738\leq \left(k -1,596.67\right)/518.738\right)  \\
&\approx&\exp(-e^{-\mu})=0.90, 
\end{eqnarray*} 
where \, $\mu=\left(k -1596.67\right)/518.738.$ Similarly, with probability $0.90$ one needs at least $\textbf{2,765}$ trials to collect all $100$ different types of coupons. 
\smallskip

Finally, suppose that all the $N$ different coupons are uniformly distributed. Hence, (\ref{N8c}) yields that with probability $0.90,$ at least $\textbf{686}$ trials are needed.

\end{document}